\documentclass[12pt]{article}
\usepackage{latexsym}
\usepackage{amssymb}
\usepackage{amsmath}
\usepackage{longtable}

\def\proof{\medbreak\noindent{\bf Proof}}

\newtheorem{theorem}{Theorem}

\def\{{\lbrace}
\def\}{\rbrace}

\def\cl{{\cal C}\!\ell}

\def\R{{\Bbb R}}

\def\C{{\Bbb C}}

\def\H{{\Bbb H}}

\def\K{{\cal K}}

\def\be{\begin{equation}}
\def\ee{\end{equation}}

\def\Mat{{\rm Mat}}

\def\Even{{\rm Even}}
\def\Odd{{\rm Odd}}

\def\be{\begin{equation}}
\def\ee{\end{equation}}

\newcommand{\fin}{\hbox{$\blacksquare$}}

\begin{document}

\title{Classification of extended Clifford algebras}


\author{N.~G.~Marchuk}

\maketitle

\thanks{nmarchuk@mi.ras.ru\\Steklov Mathematical Institute of RAS,\\ Gubkina st.8, Moscow, 119991, Russia}

\begin{abstract}
Considering  tensor products of  special commutative algebras and  general real Clifford algebras, we arrive at extended Clifford algebras.
 We have found that there are five types of extended Clifford algebras. The class of extended Clifford algebras is closed with respect to the tensor product.
\end{abstract}

MSC class:	15A66
\medskip

Key words: Clifford algebra, Quasi Clifford algebras, extended Clifford algebra, tensor product, Cartan-Bott periodicity, complexification.

\bigskip

Considering  tensor products of  special commutative algebras and  general real Clifford algebras, we arrive at extended Clifford algebras $\cl(r,s|p,q)$. This class of algebras is a subclass of the class of so-called Quasi Clifford algebras \cite{Gas}. A special case of extended Clifford algebras (in our notations $\cl(m,0|0,n)$) is used in \cite{Rajan} for coding theory.  In particular, the class of extended Clifford algebras contains the class of all Clifford algebras and the classes of complexified Clifford algebras $\C\otimes\cl(p,q)$ and hyperbolically complexified Clifford algebras $(\R\oplus\R)\otimes\cl(p,q)$ \cite{Porteous},\cite{MarShirROMP}. The class of extended Clifford algebras is closed with respect to the operation of tensor product. We develop some results of our previous paper  \cite{jreduc1} where we study tensor products of Clifford algebras using the Jacobson's approach \cite{Jac}.  A classification of extended Clifford algebras is discussed. We will see that there are five types of extended Clifford algebras and we will find corresponding conditions for integers $r,s,p,q$.


\section{Clifford algebras $\cl(p,q)$, commutative algebras $\K(r,s)$, and extended Clifford algebras $\cl(r,s|p,q)$}

\noindent{\bf Tensor products of algebras}
Let ${\cal A}$ be a real algebra, which are the real vector space $V$ with a product of elements $A_1,A_2\to
A_1\cdot A_2$ that satisfies the distributive condition. And let ${\cal B}$ be a real algebra, which are the real vector space $W$ with a product of elements $B_1,B_2\to B_1*B_2$ that satisfies the distributive condition.
A {\em tensor product} ${\cal A}\otimes{\cal B}$ of algebras ${\cal A}$ and ${\cal B}$ is the algebra, which is the tensor product $V\otimes W$ of vector spaces $V$ and $W$ with the operation of multiplication
\begin{equation}
(A_1\otimes B_1)\circ (A_2\otimes B_2) = (A_1\cdot
A_2)\otimes(B_1*B_2), \quad \forall A_1,A_2\in{\cal A},\forall
B_1,B_2\in{\cal B}.\label{tensor:product:def}
\end{equation}
\medskip

\noindent{\bf Isomorphic real algebras.}
Real algebras  ${\cal A}$ and ${\cal B}$ of the same dimension are called {\em isomorphic} if there exists a one-to-one correspondence (bijection) $T : {\cal A} \to {\cal B}$ such that
\begin{eqnarray}
T(\lambda_1 A_1 + \lambda_2 A_2) &=& \lambda_1 T(A_1)+\lambda_2 T(A_2),\quad \forall A_1,A_2\in{\cal A},\quad \forall\lambda_1,\lambda_2\in\R\nonumber\\
T(A_1\cdot A_2) &=& T(A_1)*T(A_2),\quad\forall A_1,A_2\in{\cal A}.\label{isomorpism}
\end{eqnarray}



\noindent{\bf Clifford algebras $\cl(p,q)$.}
Consider a real Clifford algebra $\cl(p,q)$  of signature $(p,q)$ with $n$ generators $n=p+q$, $n\geq1$ (the construction of Clifford algebra is discussed in details in \cite{Lounesto}).
Let $e$ be the identity element and let $e^a$, $a=1,\ldots,n$ be generators of Clifford algebra $\cl(p,q)$,\label{Cl:def}
\begin{equation}
e^a e^b + \omega e^b e^a=0,\quad\forall a,b=1,\ldots,n,\quad a\neq b,\label{com:anticom1}
\end{equation}
\begin{equation}
(e^a)^2=\left\lbrace
\begin{array}{rl}
e, & \mbox{\rm for $a=1,\ldots,p$,}\\
-e, & \mbox{\rm for $a=p+1,\ldots,n$,}
\end{array}
\right.\label{ea:square1}
\end{equation}
where $\omega=1$. The dimension of the Clifford algebra $\cl(p,q)$ (as a vector space) is equal to $2^n$. Let us define
$\cl(0,0) :=\R$.
\medskip

\noindent{\bf Commutative algebras $\K(p,q)$.}
If we take the above conditions (\ref{com:anticom1}) and (\ref{ea:square1}) in the definition of the Clifford algebra $\cl(p,q)$ and replace the identity  $\omega=1$ by the identity  $\omega=-1$, then we arrive at the commutative algebra $\K(p,q)$ of signature $(p,q)$ with $n$ generators $n=p+q$, $n\geq1$.
The dimension of the algebra $\K(p,q)$ (as a vector space) is equal to $2^n$.
 Also define
$\K(0,0):=\R$.

\medskip

\noindent{\bf Tensor products of Clifford algebras.} Let $p_1,\ldots,p_k$,$q_1,\ldots,q_k$ be nonnegative integer numbers and  $n_j=p_j+q_j$, $j=1,\ldots k$ be natural numbers. And let $m$ be the number of odd numbers in the set $n_1,\ldots,n_k$. Consider the tensor product of Clifford algebras
\begin{equation}
\cl(p_1,q_1)\otimes\ldots\otimes\cl(p_k,q_k).\label{tpca}
\end{equation}
The following propositions are consequences of Theorem 4 and Theorem 5 from \cite{jreduc1}:
\begin{enumerate}
\item If $m=0$ or $m=1$, then there exists a nonnegative integer numbers $p,q$ such that $p+q=n=n_1+\ldots+n_k$ and the tensor product (\ref{tpca}) is isomorphic to the Clifford algebra $\cl(p,q)$;
    \item If $m\geq 2$, then there exists a nonnegative integer numbers $r,s,p,q$ such that $r+s=m$, $p+q=n=n_1+\ldots+n_k-m$ and the tensor product (\ref{tpca}) is isomorphic to the tensor product
        $$
        \K(r,s)\otimes\cl(p,q),
        $$
        which can be considered as the Clifford algebra $\cl(p,q)$ over the commutative algebra (ring) $\K(r,s)$.
\end{enumerate}

\medskip

\noindent{\bf Extended Clifford algebras.} The above two propositions lead us to consider a class of  {\em extended Clifford algebras}
\begin{equation}
\cl(r,s|p,q) := \K(r,s)\otimes\cl(p,q),\label{ex:cl}
\end{equation}
where $r,s,p,q$ -- nonnegative integer numbers. Any algebra from this class is generated by a set of generators, which is the union of two sets -- a set of commuting generators and a set of anticommuting generators. Generators from different sets are commute. As a result we arrive at a class of associative unital algebras dependent on four nonnegative integer numbers. The dimension of this algebra (as a vector space) is equal to  $2^{r+s+p+q}$. In particular, the class of extended Clifford algebras contains the class of Clifford algebras $\cl(p,q)$, the class of complexified Clifford algebras $\C\otimes\cl(p,q)$, and the class of hyperbolically complexified Clifford algebras $(\R\oplus\R)\otimes\cl(p,q)$.

We see that the class of extended Clifford algebras
 is closed w.r.t. the tensor product, i.e., for any tensor product
 \begin{equation}
 \cl(r_1,s_1|p_1,q_1)\otimes\ldots\otimes\cl(r_k,s_k|p_k,q_k)\label{tpea}
 \end{equation}
there exist four nonnegative integer numbers $r,s,p,q$ such that $r+s+p+q=\sum_{j=1}^k(r_j+s_j+p_j+q_j)$ and the tensor product (\ref{tpea}) is isomorphic to the extended Clifford algebra $\cl(r,s|p,q)$.

\section{Classification of extended Clifford algebras}
Let $M\geq1,N\geq0$ be  integer numbers. Consider five algebras each of which is a tensor product of Clifford algebras
with one or two generators
\begin{enumerate}
\item[(I)] $ $ $\cl(1,1)^N$;
\item[(II)]$ $  $\cl(0,2)\otimes\cl(1,1)^{N-1}$, $N\geq1$;
\item[(III)] $ $ $\cl(0,1)^M\otimes\cl(1,1)^N$;
\item[(IV)] $ $ $\cl(1,0)^M\otimes\cl(1,1)^N$;
\item[(V)]$ $  $\cl(1,0)^M\otimes\cl(0,2)\otimes\cl(1,1)^{N-1}$, $N\geq1$,
\end{enumerate}
where by $\cl(p,q)^M$ we denote the tensor product of $M$ pieces of Clifford algebra $\cl(p,q)$ and
$$
\cl(p,q)^0 := \cl(0,0)\simeq\R.
$$
 Algebras  I-V depend on one ($N$) or two ($M,N$) nonnegative integer numbers. In the sequel we say about algebras of type I,II,III,IV,V.

\begin{theorem}
Suppose that  $r,s,p,q$ are nonnegative integer numbers and
\begin{equation}
m=r+s,\quad n=p+q,\quad M=m+\sigma(n),\quad N=[\frac{n}{2}],\quad t= (p-q)\,{\rm mod}\,8,\label{mnMN}
\end{equation}
where $[\frac{n}{2}]$ is the integer part of  $n/2$ and
\begin{equation}
\sigma(n):=n-2[\frac{n}{2}]=\left\lbrace
\begin{array}{rl}
0, & \mbox{\rm for even $n$,}\\
1, & \mbox{\rm for odd $n$.}
\end{array}
\right.\nonumber
\end{equation}
Then the extended Clifford algebra $\cl(r,s|p,q)$ is isomorphic to one of tensor products of type {\rm I-V} according to the table\footnote{For types III,IV,V we have two lines in the table that means there are two cases. For example, we have type IV for $r\geq0$, $s=0$, $t=1$ or for $r\geq1$, $s=0$, $t=0,2$.}
$$
\begin{tabular}{|c||c|c|c|}
\hline
{\rm type} & $r$ & $s$ & $t$    \\ \hline \hline
{\rm I}    & $0$ & $0$ & $0,2$  \\ \hline
{\rm II}   & $0$ & $0$ & $4,6$  \\ \hline
{\rm III}  & $\geq0$ & $\geq1$ & $0,1,2,3,4,5,6,7$ \\
     & $\geq0$ & $0$ & $3,7$ \\ \hline
{\rm IV}   & $\geq0$ & $0$ & $1$  \\
     & $\geq1$ & $0$ & $0,2$ \\ \hline
{\rm V}    & $\geq0$ & $0$ & $5$ \\
     & $\geq1$ & $0$ & $4,6$\\
     \hline
\end{tabular}
$$
\end{theorem}

\proof. For Clifford algebras $\cl(p,q)$ we have the following Cartan's classification \cite{jreduc1}:
\begin{equation}
\cl(p,q)\simeq\left\lbrace
\begin{array}{ll}
\cl(1,1)^{\frac{n}{2}}, & \parbox{.5\linewidth}{ for $p-q=0, 2\!\!\mod 8$;}\\
\cl(0,2)\otimes\cl(1,1)^{\frac{n-2}{2}}, & \parbox{.5\linewidth}{ for $p-q=4, 6\!\!\mod 8$;}\\
\cl(0,1)\otimes\cl(1,1)^{\frac{n-1}{2}}, & \parbox{.5\linewidth}{ for $p-q=3, 7\!\!\mod 8$;}\\
\cl(1,0)\otimes\cl(1,1)^{\frac{n-1}{2}}, & \parbox{.5\linewidth}{ for $p-q=1\!\!\mod 8$;}\\
\cl(1,0)\otimes\cl(0,2)\otimes\cl(1,1)^{\frac{n-3}{2}}, & \parbox{.5\linewidth}{ for $p-q=5\!\!\mod 8$.}
\end{array}
\right.\label{my:formula1}
\end{equation}
If $r=s=0$ and $t$ is even, then
$$
\cl(r,s|p,q)\simeq\cl(p,q)\simeq\left\lbrace
\begin{array}{ll}
\cl(1,1)^{\frac{n}{2}}, & \mbox{\rm for $t=0,2$,}\\
\cl(0,2)\otimes\cl(1,1)^{\frac{n}{2}-1}, & \mbox{\rm for  $t=4,6$.}
\end{array}
\right.\nonumber
$$
That means $\cl(0,0|p,q)$ are of the type I for $t=0,2$ and of the type II for $t=4,6$.

If $s=0$, $r=m\geq1$, then $\K(r,s)\simeq\cl(1,0)^m$ and from (\ref{my:formula1}) we get
$$
\cl(r,s|p,q)\simeq\cl(1,0)^m\otimes\cl(p,q)
$$
and
\begin{equation}
\cl(r,s|p,q)
\simeq\left\lbrace
\begin{array}{ll}
\cl(1,0)^m\otimes\cl(1,1)^{\frac{n}{2}}, & \parbox{.5\linewidth}{ for $t=0, 2$;}\\
\cl(1,0)^m\otimes\cl(0,2)\otimes\cl(1,1)^{\frac{n-2}{2}}, & \parbox{.5\linewidth}{ for $t=4, 6$;}\\
\cl(0,1)^{m+1}\otimes\cl(1,1)^{\frac{n-1}{2}}, & \parbox{.5\linewidth}{ for $t=3, 7$;}\\
\cl(1,0)^{m+1}\otimes\cl(1,1)^{\frac{n-1}{2}}, & \parbox{.5\linewidth}{ for $t=1$;}\\
\cl(1,0)^{m+1}\otimes\cl(0,2)\otimes\cl(1,1)^{\frac{n-3}{2}}, & \parbox{.5\linewidth}{ for $t=5$.}
\end{array}
\right.\label{cl:rspq}
\end{equation}
That means each line at the right hand part of (\ref{cl:rspq}) is of the type IV,V,III,IV,V respectively. For $t=3,7$ we use the isomorphism $\cl(1,0)^m\otimes\cl(0,1)\simeq\cl(0,1)^{m+1}$.

If $s\geq1$, then $\K(r,s)\simeq\cl(0,1)^m$. Using the isomorphism
$$
\cl(0,1)\otimes\cl(0,2)\simeq\cl(0,1)\otimes\cl(1,1),
$$
we see that
$$
\cl(r,s|p,q)\simeq\cl(0,1)^M\otimes\cl(1,1)^N
$$
are of the type III for all $t=0,1,2,3,4,5,6,7$.
This completes the proof.
\medskip

\noindent{\bf Consequence 1}. Suppose that we have two sets of four nonnegative integer numbers
 $r,s,p,q$ and $\acute r,\acute s,\acute p,\acute q$ such that $r+s+p+q= \acute r+\acute s+\acute p+\acute q$. Denote
$$
m=r+s,\quad n=p+q,\quad M=m+\sigma(n),\quad N=[\frac{n}{2}],
$$
$$
\acute m=\acute r+\acute s,\quad \acute n=\acute p+\acute q,\quad \acute M=\acute m+\sigma(\acute n),\quad \acute N=[\frac{\acute n}{2}].
$$
Extended Clifford algebras
 $\cl(r,s|p,q)$ and $\cl(\acute r,\acute s|\acute p,\acute q)$ are isomorphic iff $\acute M=M$, $\acute N=N$ and both
 extended Clifford algebras belong to the same type (I-V).

\medskip
\noindent{\bf Consequence 2}. For any fixed natural $n$ all $n+1$ complexified Clifford algebras $\C\otimes\cl(p,q)$, $p+q=n$ are isomorphic\footnote{This proposition is evident for complex Clifford algebras because we can take new generators $\acute e_k=i e_k$. But for complexified Clifford algebras this proposition is not so evident (at least for me).}.
\medskip

\noindent{\bf Example.} Let us consider the following question\footnote{The answer to this question can be found using the theory of semisimple algebras -- you need to do some efforts for this.}: are there isomorphic algebras among the extended Clifford algebras
\begin{equation}
\cl(3,0|7,15),\quad\cl(4,0|3,18),\quad\cl(5,0|11,9) ?\label{3alg}
\end{equation}
We can immediately  get the answer on the question. Let $m_i,n_i,M_i,N_i,t_i$, $i=1,2,3$ be the numbers defined for algebras (\ref{3alg}) by formulas (\ref{mnMN})
\begin{eqnarray*}
&& m_1=3,\ n_1=22,\ M_1=3,\ N_1=11,\ t_1=0,\\
&& m_2=4,\ n_2=21,\ M_2=5,\ N_2=10,\ t_2=1,\\
&& m_3=5,\ n_3=20,\ M_3=5,\ N_3=10,\ t_3=2.
\end{eqnarray*}
 We see that $M_2=M_3,\ N_2=N_3$ and both extended Clifford algebras $\cl(4,0|3,18), \cl(5,0|11,9)$ belong to the same type IV. According to consequence 1 these algebras are isomorphic.

\medskip

\noindent{\bf Acknowledgments.}
I'm grateful to Prof P.Leopardi for references to important papers \cite{Gas},\cite{Rajan}.
  
  This work was supported by the grant of the Russian
Science Foundation RSF 14-50-00005.


\end{document}